# APPROXIMATION BY NONLINEAR FOURIER BASIS IN GENERALIZED HÖLDER SPACES

Hatice ASLAN[*1], Ali GUVEN[2]

[*1] *Firat University, Department of Mathematics,23119 Elazig/Turkey; ORCID: 0000-0002-3486-4179*
[2] *Balikesir University, Department of Mathematics,10145, Balikesir/Turkey; ORCID: 0000-0001-8878-250X*

**ABSTRACT**

In this paper, the value of deviation of a function $f$ from its $n$th generalized de la Vallèe-Poussin mean $V_n^a(\lambda, f)$ with respect to the nonlinear trigonometric system is estimated for the classes of $2\pi$-periodic functions in the uniform norm $\|.\|_{C(\mathbb{T})}$ and in the generalized Hölder norm $\|.\|_{\omega_\beta}$ where $f \in H^{\omega_\alpha}(\mathbb{T})$ and $0 \leq \beta < \alpha \leq 1$.

## 1. INTRODUCTION

Trigonometric (or equivalently exponential) Fourier series are most common tools for studying approximation properties of $2\pi$-periodic functions on the real line. There are many methods such as summation of partial sums of Fourier series which are used in approximation theory (Cesàro, Abel-Poisson means, de la Vallèe-Poussin, etc.). These methods are used by several authors to study approximation properties of functions in $C(\mathbb{T})$ (the space of $2\pi$-periodic continuous functions on $\mathbb{R}$), $H^\alpha(\mathbb{T})$ and $H^{\omega_\alpha}(\mathbb{T})$ (the Hölder class and generalized Hölder class of $2\pi$-periodic functions, where $0 < \alpha \leq 1$ and $\omega_\alpha \in \mathcal{M}_\alpha$). Results of these studies can be found in the monographs [1-4], and in the survey [5]. Furthermore some kinds of results "Korovkin-type theorems" discovered such a property for the functions $1, \cos$ and $\sin$ in $C(\mathbb{T})$ by Korovkin in 1953. These theorems of Korovkin are actually equivalent to the trigonometric version of the classical Weierstrass approximation theorem. These theorems exhibit a variety of test subsets of functions which guarantee that the approximation (or the convergence) property holds on the whole space provided it holds on them. After his discovery, Korovkin-type approximation theory have been extending by several mathematicians [6-7]. Also, for some special classes, namely Hölder classes of continuous $2\pi$-periodic functions, there are several approximation results. For example, S. Prössdorf studied the degree of approximation of Cesàro means of Fourier series of functions in Hölder classes [8] and Z. Stypinski used generalized de la Vallèe-Poussin means of Fourier series and extended results of S. Prössdorf [9]. Later, Leindler defined more general classes than Hölder classes (the generalized Hölder classes) and investigated approximation properties of generalized de la Vallèe-Poussin means of series [1]. Results of these studies and more results on Hölder approximation can be found in [10].

A family of nonlinear Fourier bases $e^{ik\theta_a(x)}, k \in \mathbb{Z}$, as the extension of the classical Fourier basis, defined by the nontangential boundary value of the Möbius transformation and applied to signal processing [11-15].

In [16], the author established the Jackson's and Bernstein's theorems for approximation of functions in $L^p(\mathbb{T}), 1 \leq p \leq \infty$, by the nonlinear Fourier basis. In [17], the order of convergence of classical and generalized de la Vallèe-Poussin means of of Fourier series by nonlinear basis of continuous functions is investigated in uniform and Hölder norms.

The aim of this article is to obtain estimates for the approximation order of generalized de la Vallée Poussin means of series with nonlinear Fourier basis in uniform and generalized Hölder norms.

## 2. SERIES BY NONLINEAR FOURIER BASIS

In this section, we shall give the definition and basic properties of series by nonlinear Fourier basis. We also give the definition of the generalized Hölder classes.

Let $\mathbb{D} = \{z \in \mathbb{C}: |z| < 1\}$ and $a \in \mathbb{D}$. We consider the Möbius transformation

$$\tau_a(z) = \frac{z-a}{1-\bar{a}z},$$

which is a conformal automorphism of $\mathbb{D}$. The nonlinear phase function $\theta_a$ is defined through the relation

$$e^{ik\theta_a(t)} := \tau_a(e^{it}) = \frac{e^{it}-a}{1-\bar{a}e^{it}},$$

where $\tau_a(e^{it})$ stands for radial boundary value of $\tau_a$. It is easy to see that $\theta_a(t+2\pi) = \theta_a(t) + 2\pi$, and if we set $a = |a|e^{it_a}$ then

$$\theta'_a(t) = p_a(t) := \frac{1-|a|^2}{1-2|a|\cos(t-t_a)}$$

which satisfies

$$\frac{1-|a|}{1+|a|} \leq p_a(t) \leq \frac{1+|a|}{1-|a|}.$$

Let $\mathbb{T} = \mathbb{R}/2\pi\mathbb{Z}$. The space $L_a^2(\mathbb{T})$ consists of measurable functions $f: \mathbb{T} \to \mathbb{C}$ such that

$$\frac{1}{2\pi}\int_{\mathbb{T}} |f(t)|^2 p_a(t) dt < \infty$$

becomes a Hilbert space with respect to the inner product

$$<f,g>_a := \frac{1}{2\pi}\int_{\mathbb{T}} f(x)\overline{g(x)} p_a(x) dx$$

and the set

$$\{e^{ik\theta_a(x)}: k \in \mathbb{Z}\}$$

is an orthonormal basis for $L_a^2(\mathbb{T})$ (see [11] and [12]). In the case $a = 0$ we obtain the classical Fourier basis $\{e^{ikx}: k \in \mathbb{Z}\}$ for the space $L_0^2(\mathbb{T}) = L^2(\mathbb{T})$.

Fourier series of a function $f \in L^2(\mathbb{T})$ with respect to the orthonormal basis $\{e^{ik\theta_a(x)}: k \in \mathbb{Z}\}$ become

$$f(x) \sim \sum_{-\infty}^{\infty} c_k^a(f) e^{ik\theta_a(x)}, \tag{2.1}$$

where

$$c_k^a(f) := \frac{1}{2\pi}\int_{\mathbb{T}} f(x) e^{ik\theta_a(x)} p_a(x) dx, k \in \mathbb{Z}.$$

If we denote the $n$th partial sum of this series by $S_n^a(f)(x)$ then we have

$$S_n^a(f)(x) = \frac{1}{2\pi}\int_{\mathbb{T}} f(t) D_n(\theta_a(x) - \theta_a(t)) p_a(t) dt$$

$$= \frac{1}{2\pi}\int_{\mathbb{T}} f(t) F(\theta_a(x) + t) D_n(t) dt$$

where $F := f \circ \theta_a^{-1}$ and $D_n$ is the Dirichlet kernel of order $n$.

For each natural number $n$, let $\tau_n^a$ be the set of all nonlinear trigonometric polynomials of degree at most $n$, that is

$$\tau_n^a := \text{span}\{e^{ik\theta_a(x)} : |k| \leq n\},$$

and let $E_n^a(f)$ be the approximation error of $f \in C(\mathbb{T})$ by elements of $\tau_n^a$ i. e.

$$E_n^a(f) := \inf_{T \in \Pi_n^a} \|f - T\|_\infty.$$

Let $\lambda = \{\lambda_n\}$ be a sequence of integers such that $\lambda_1 = 1$ and $0 \leq \lambda_{n+1} - \lambda_n \leq 1$. The sequence of generalized de la Vallée Pousin means of the series (2.1) is defined by

$$V_n^a(\lambda, f) = \frac{1}{\lambda_n} \sum_{k=n-\lambda_n}^{n-1} S_k^a(f).$$

In special case $\lambda_n = 1$ $(n = 1,2,...)$ $V_n^a(\lambda, f)$ become $S_{n-1}^a(f)$ and in the case $\lambda_n = n$ $(n = 1,2,...)$ have the Fejér means

$$V_n^a(\lambda, f) = \sigma_{n-1}^a(f) := \frac{1}{n} \sum_{k=0}^{n-1} S_k^a(f).$$

In the special case, $\lambda_n = n$ the means $V_n^a(\lambda, f)$ coincide with the Cesàro $(C, 1)$ means of (2.1).

We refer to [11] and [16] for more detailed information on Fourier series by nonlinear basis.

## 3. APPROXIMATION IN THE GENERALIZED HÖLDER NORM

We denote by $C(\mathbb{T})$ the Banach space of continuous functions $f: \mathbb{R} \to \mathbb{C}$ equipped with the norm

$$\|f\|_\infty := \sup_{x \in \mathbb{T}} |f(x)|.$$

The modulus of continuity of $f \in C(\mathbb{T})$ is defined by

$$\omega(f, t) := \sup_{\substack{x,y \in \mathbb{T} \\ |x-y| \leq t}} |f(x) - f(y)|$$

for $t > 0$.

For any modulus of continuity $\omega$, we define the generalized Hölder class $H^\omega(\mathbb{T})$ as the set of functions $f \in C(\mathbb{T})$ for which

$$A^\omega(f) := \sup_{t \neq s} \frac{|f(t) - f(s)|}{\omega(\|t - s\|)} < \infty,$$

and the norm on $H^\omega(\mathbb{T})$ as

$$\|f\|_\omega := \|f\|_\infty + A^\omega(f).$$

If $\omega(\delta) = \delta^\alpha, 0 < \alpha \leq 1$, then we write $H^\alpha(\mathbb{T})$ instead of $H^\omega(\mathbb{T})$ and $\|f\|_\alpha$ instead of $\|f\|_\omega$.

In [1], L. Leindler introduced a certain class of moduli of continuity for $0 \leq \alpha \leq 1$, let $\mathcal{M}_\alpha$ denote the class of moduli of continuity $\omega_\alpha$ having the following properties:

(i) for any $\alpha' > \alpha$, there exists a natural number $\mu = \mu(\alpha')$ such that

$$2^{\mu\alpha'} \omega_\alpha(2^{-n-\mu}) > 2\omega_\alpha(2^{-n}), (n = 1,2,...),$$

($ii$) for every natural number $v$, there exists a natural number $N(v)$ such that

$$2^{v\alpha'}\omega_\alpha(2^{-n-v}) > 2\omega_\alpha(2^{-n}), (n > N(v)).$$

It is clear that $\omega(\delta) = \delta^\alpha \in \mathcal{M}_\alpha$ but $\omega_\alpha(\delta)$ is an extension of $\omega(\delta) = \delta^\alpha$. Consequently, in general, $H^{\omega_\alpha}(\mathbb{T})$ is larger than $H^\alpha(\mathbb{T})$.

It is known that if $0 \leq \beta < \alpha \leq 1$, $\omega_\beta \in \mathcal{M}_\beta$ and $\omega_\alpha \in \mathcal{M}_\alpha$ then the function $\omega_\alpha/\omega_\beta$ is non-decreasing [18].

We shall use the notation $A \lesssim B$ at inequalities if there exists an absolute constant $c > 0$, such that $A \leq cB$ holds for quantities $A$ and $B$.

## 4. AUXILIARY RESULTS

**Lemma 1** [1] If $0 \leq \beta < \alpha \leq 1$, $\omega_\beta \in \mathcal{M}_\beta$, $\omega_\alpha \in \mathcal{M}_\alpha$ and $f \in H^{\omega_\alpha}(\mathbb{T})$, then

$$\sum_{k=1}^n \frac{\omega_\beta(2^{-k})}{\omega_\alpha(2^{-k})} \leq K_{\alpha,\beta} \frac{\omega_\beta(2^{-n})}{\omega_\alpha(2^{-n})} \tag{4.1}$$

and

$$\sum_{k=n}^\infty \frac{\omega_\beta(2^{-k})}{\omega_\alpha(2^{-k})} \leq K_{\alpha,\beta} \frac{\omega_\beta(2^{-n})}{\omega_\alpha(2^{-n})} \tag{4.2}$$

hold where $K_{\alpha,\beta}$ a positive constant independent of $n$.

**Lemma 2** [1] For any nonnegative sequence $a_n$, the inequality

$$\sum_{n=1}^m a_n \leq K a_m, m = 1,2,\ldots; K > 0$$

holds if and only if there exist a positive number $c$ and a natural number $\mu$ such that for any $n$,

$$a_{n+1} \leq c a_n$$

and

$$a_{n+1} \geq 2 a_n$$

are valid.

Let $\varphi$ be an increasing positive function on $(0,\infty)$. The $\varphi$-norm of a function $f \in C(\mathbb{T})$ is defined by

$$\|f\|_\varphi := \|f\|_\infty + \sup_{x \neq y} \frac{|f(x) - f(y)|}{\varphi(|x-y|)} = \|f\|_\infty + \sup_{\delta > 0} \frac{\|f - f(. + \delta)\|_\infty}{\varphi(\delta)}.$$

It is clear that, special case $\varphi(\delta) = \delta^\alpha (0 < \alpha \leq 1)$, we have $\|f\|_\varphi = \|f\|_\alpha$.

The following important result was obtained in [19].

**Lemma 3** Let $\{A_n\}$ be a sequence of linear convolution operators from $C(\mathbb{T})$ into $C(\mathbb{T})$ and let $\varphi$ be an increasing positive function on $(0,\infty)$. Then

$$\|A_n(f) - f\|_\varphi := \|A_n(f) - f\|_\infty \left(1 + \frac{2}{\varphi\left(\frac{1}{n}\right)}\right) + \sup_{0 < \delta \leq \frac{1}{n}} \frac{2\omega(f, \delta)}{\varphi(\delta)} (1 + \|A_n\|) \tag{4.3}$$

holds for every $f \in C(\mathbb{T})$, where $\|A_n\|$ is the operator norm of $A_n$.

**Lemma 4** [16]  Let $f \in C(\mathbb{T})$.  Then we have

$$\frac{1-|a|}{2} \omega(f, \delta)_\infty \leq \omega(f \circ \theta_a^{-1}, \delta)_\infty \leq \frac{2}{1-|a|} \omega(f, \delta)_\infty.$$

## 5. MAIN RESULTS

Our main results are following.

Let $\lambda = \{\lambda_n\}$ be a sequence of integers such that $\lambda_1 = 1$ and $0 \leq \lambda_{n+1} - \lambda_n \leq 1$. The sequence of generalized de la Vallée Pousin means of the series (2.1) is defined by

$$V_n^a(\lambda, f) = \frac{1}{\lambda_n} \sum_{k=n-\lambda_n}^{n-1} S_k^a(f).$$

In special case $\lambda_n = 1$ $(n = 1, 2, ...)$ $V_n^a(\lambda, f)$ become $S_{n-1}^a(f)$ and in the case $\lambda_n = n$ $(n = 1, 2, ...)$ have the Fejér means

**Theorem 1** [17] Let $\lambda = \{\lambda_n\}$ be a sequence of integers such that $\lambda_1 = 1$ and $0 \leq \lambda_{n+1} - \lambda_n \leq 1$. If $f \in C(\mathbb{T})$ satisfies $|f(x)| < M$, then the estimate

$$\|V_n^a(\lambda, f)\|_\infty \leq M \left(3 + \log \frac{2n - \lambda_n}{\lambda_n}\right) \tag{5.1}$$

holds for every natural number $n$.

**Theorem 2** [17] $f \in C(\mathbb{T})$ the degree of approximation by the sequence $\{V_n^a(\lambda, f)\}$ of generalized de la Vallée Poussin means is estimated as

$$\|f - V_n^a(\lambda, f)\|_\infty \lesssim \left(3 + \log \frac{2n - \lambda_n}{\lambda_n}\right) E_{n-\lambda_n}^a(f).$$

We obtained the following estimation for the deviation of $V_n^a(\lambda, f)$ from $f \in H^{\omega_\alpha}(\mathbb{T})$ in the uniform norm.

**Theorem 3**  Let $0 \leq \beta < \alpha \leq 1$, $\omega_\beta \in \mathcal{M}_\beta$, $\omega_\alpha \in \mathcal{M}_\alpha$ and $f \in H^{\omega_\alpha}(\mathbb{T})$. Then

$$\|f - V_n^a(\lambda, f)\|_\infty \lesssim \begin{cases} \left(\frac{12}{1-|a|}\right) \omega_\alpha(1/\lambda_n), & \alpha < 1 \text{ or } \beta > 0 \\ \left(\frac{12}{1-|a|}\right) \omega_\alpha(1/\lambda_n)(1 + \log \lambda_n), & \alpha = 1 \text{ and } \beta = 0 \end{cases} \tag{5.2}$$

The estimation of $f - V_n^a(\lambda, f)$ in the generalized Hölder norm is obtained in the following theorem.

**Theorem 4**  Let $0 \leq \beta < \alpha \leq 1$, $\omega_\beta \in \mathcal{M}_\beta$, $\omega_\alpha \in \mathcal{M}_\alpha$ and $f \in H^{\omega_\alpha}(\mathbb{T})$. Then the estimate

$$\|f - V_n^a(\lambda, f)\|_{\omega_\beta} \lesssim \begin{cases} \left(\frac{12}{1-|a|}\right) \frac{\omega_\beta\left(\frac{1}{n}\right)}{\omega_\alpha\left(\frac{1}{n}\right)} \omega_\alpha(1/\lambda_n) \log \frac{2n}{\lambda_n}, & \alpha < 1 \text{ or } \beta > 0 \\ \left(\frac{12}{1-|a|}\right) \frac{\omega_\beta\left(\frac{1}{n}\right)}{\omega_\alpha\left(\frac{1}{n}\right)} \omega_\alpha(1/\lambda_n)(1 + \log \lambda_n) \log \frac{2n}{\lambda_n}, & \alpha = 1 \text{ and } \beta = 0 \end{cases} \tag{5.3}$$

*holds.*

**Corollary 1** *Let If* $0 \leq \beta < \alpha \leq 1$ *and* $f \in H^\alpha(\mathbb{T})$. *Then*

$$\|f - V_n^a(\lambda, f)\|_\beta \lesssim \begin{cases} \left(\dfrac{1+|a|}{1-|a|}\right)^\alpha n^{\beta-\alpha} \dfrac{1}{\lambda_n^\alpha} \log \dfrac{2n}{\lambda_n}, & \alpha < 1 \\ \left(\dfrac{1+|a|}{1-|a|}\right) n^{\beta-1} \left(\dfrac{1+\log \lambda_n}{\lambda_n}\right) \log \dfrac{2n}{\lambda_n}, & \alpha = 1. \end{cases}$$

## 6. PROOFS OF MAIN RESULTS

*Proof of Theorem 3.* Let

$$\phi_x^a(2t) = F(\theta_a(x) + 2t) - F(\theta_a(x) - 2t) - 2f(x)$$

where $F := f \circ \theta_a^{-1}$. A standart computation gives that

$$V_n^a(\lambda, f) - f(x) = \frac{1}{\lambda_n \pi} \int_0^{\frac{\pi}{2}} \phi_x^a(2t) K_n(t) dt$$

where

$$K_n(t) = \frac{\sin \lambda_n t \sin(2n - \lambda_n)t}{\sin^2 t}.$$

Hence,

$$|V_n^a(\lambda, f) - f(x)| \lesssim \frac{1}{\lambda_n \pi} \int_0^{\frac{\pi}{2}} |\phi_x^a(2t)| K_n(t) dt.$$

Therefore, since $f \in H^\alpha(\mathbb{T})$ from Lemma 4 it is clear that

$$|\phi_x^a(2t)| \lesssim \left(\frac{12}{1-|a|}\right) \omega_\alpha(t). \tag{6.1}$$

Let us split the integral $I$ into three parts and by using Lemma 5:

$$I := \int_0^{\frac{\pi}{2}} = \int_0^{\frac{1}{2n-\lambda_n}} + \int_{\frac{1}{2n-\lambda_n}}^{\frac{1}{\lambda_n}} + \int_{\frac{1}{\lambda_n}}^{\frac{\pi}{2}} = I_1 + I_2 + I_3$$

These integrals by (4.1), (4.2) and (6.1) can be estimated by elementary methods:

$$\frac{1}{\lambda_n} I_1 \leq \frac{1}{\lambda_n} \int_0^{\frac{1}{2n-\lambda_n}} \lambda_n |\phi_x^a(2t)| dt \lesssim \left(\frac{12}{1-|a|}\right) \lambda_n \omega_\beta\left(\frac{1}{2n-\lambda_n}\right) \int_0^{\frac{1}{2n-\lambda_n}} \frac{\omega_\alpha(t)}{\omega_\beta(t)} dt$$

$$\lesssim \left(\frac{12}{1-|a|}\right) \omega_\beta(1/2n - \lambda_n) \frac{\omega_\alpha(1/2n - \lambda_n)}{\omega_\beta(1/2n - \lambda_n)}$$

$$\lesssim \left(\frac{12}{1-|a|}\right) \omega_\alpha(1/\lambda_n),$$

$$\frac{1}{\lambda_n} I_2 \leq \frac{1}{\lambda_n} \int_{\frac{1}{2n-\lambda_n}}^{\frac{1}{\lambda_n}} \lambda_n \frac{|\phi_x^a(t)|}{t} \lesssim \left(\frac{12}{1-|a|}\right) \lambda_n \omega_\beta\left(\frac{1}{\lambda_n}\right) \int_{\frac{1}{2n-\lambda_n}}^{\frac{1}{\lambda_n}} \frac{\omega_\alpha(t)}{t \omega_\beta(t)} dt$$

$$\lesssim \left(\frac{12}{1-|a|}\right)\omega_\beta(1/\lambda_n) \sum_{k=\lambda_n}^{2n-\lambda_n} \int_{\frac{1}{k+1}}^{\frac{1}{k}} \frac{\omega_\alpha(t)}{t\omega_\beta(t)}\,dt$$

$$\lesssim \left(\frac{12}{1-|a|}\right)\omega_\beta(1/\lambda_n) \sum_{k=\lambda_n}^{2n-\lambda_n} \frac{\omega_\alpha(1/k)}{k\omega_\beta(1/k)}$$

$$\lesssim \left(\frac{12}{1-|a|}\right)\omega_\beta(1/\lambda_n) \sum_{k=\log\lambda_n}^{\log 2n-\lambda_n} \frac{\omega_\alpha(2^{-m})}{k\omega_\beta(2^{-m})}$$

$$\lesssim \left(\frac{12}{1-|a|}\right)\omega_\beta(1/\lambda_n) \frac{\omega_\alpha(1/\lambda_n)}{\omega_\beta(1/\lambda_n)}$$

$$\lesssim \left(\frac{12}{1-|a|}\right)\omega_\alpha(1/\lambda_n),$$

and finally by using Lemma 1, we have

$$I_3 \lesssim \int_{\frac{1}{\lambda_n}}^{\frac{\pi}{2}} \frac{|\phi_x^a(t)|}{t^2}\,dt \lesssim \int_{\frac{1}{\lambda_n}}^{\frac{\pi}{2}} \left(\frac{12}{1-|a|}\right)\frac{\omega_\alpha(t)}{t^2\omega_\beta(t)}\,dt$$

$$\lesssim \left(\frac{12}{1-|a|}\right)\omega_\beta(\pi/2) \sum_{k=1}^{\lambda_n} \frac{\omega_\alpha(2^{-m})}{k\omega_\beta(2^{-m})} \lesssim \left(\frac{12}{1-|a|}\right) \sum_{k=1}^{\lambda_n} \frac{\omega_\alpha(1/k)}{\omega_\beta(1/k)}$$

$$I_3 \lesssim \left(\frac{12}{1-|a|}\right) \sum_{m=0}^{\log\lambda_n} 2^m \frac{\omega_\alpha(2^{-m})}{\omega_\beta(2^{-m})}$$

Here the last sum can be estimated easily if $\alpha < 1$; namely then

$$\sum_{m=0}^{\log\lambda_n} 2^m \frac{\omega_\alpha(2^{-m})}{\omega_\beta(2^{-m})} \lesssim \frac{1}{\omega_\beta(1/\lambda_n)} \sum_{m=0}^{\log\lambda_n} 2^m \omega_\alpha(2^{-m}) \lesssim \lambda_n \frac{\omega_\alpha(1/\lambda_n)}{\omega_\beta(1/\lambda_n)}$$

(see Lemma 1 and (5.1) with $\alpha' = 1$).

If $\alpha = 1$ and $\beta = 0$ we obtain the same upper estimation for this sum but then (5.1) holds if and only if $\alpha' > 1 (= \alpha)$. Using the monocity of the sequence $2^{m(1+\beta/2)}\omega_\alpha(2^{-m})$ we get that

$$\sum_{m=0}^{\log\lambda_n} 2^m \frac{\omega_\alpha(2^{-m})}{\omega_\beta(2^{-m})} \lesssim \lambda_n^{(1+\beta/2)}\omega_\alpha(1/\lambda_n) \sum_{m=0}^{\log\lambda_n} \frac{1}{2^{m\beta/2}\omega_\beta(2^{-m})}$$

and if we show that

$$\sum_{m=0}^{\log\lambda_n} 2^{-m\beta/2}\left(\omega_\beta(2^{-m})\right)^{-1} \lesssim \lambda_n^{-\frac{\beta}{2}}\left(\omega_\beta(2^{-m})\right)^{-1} \qquad (6.2)$$

holds, then our statement is verified.

To prove (6.2) first we show that there exists a natural number $\mu$ such that

$$2^{-m\beta/2}\left(\omega_\beta(2^{-m})\right)^{-1} > 2^{-m\beta/2}\left(\omega_\beta(2^{-m})\right)^{-1}. \qquad (6.3)$$

Since $\omega_\beta(\delta) \in \mathcal{M}_\beta$ thus, for any $\mu$, there exists an index $N(\mu)$ such that if $m > N(\mu)$ then

$$2^{m\beta}\omega_\beta(2^{-m-\mu}) < 2\omega_\beta(2^{-m});$$

and if $\mu > 4/\beta$ then hence, we get that

$$2^{m\beta}\omega_\beta(2^{-m-\mu}) < \frac{1}{2}2^{m\beta/2}\omega_\beta(2^{-m}),$$

which implies (6.2).

A standart calculation similar to the proof Lemma 2 shows that (6.3) implies (6.2).

In the case $\alpha = 1$ and $\beta = 0$ the sum investigated before does not exceed

$$K\lambda_n(1 + \log\lambda_n)\frac{\omega_1(1/\lambda_n)}{\omega_0(1/\lambda_n)}.$$

Namely, $\{2^m\omega_1(2^{-m})/\omega_0(2^{-m})\}$ is nondecreasing sequence.

Consequently, collecting the partial results, we have that

$$\frac{1}{\lambda_n}I_n \lesssim \begin{cases} \left(\frac{12}{1-|a|}\right)\omega_\alpha(1/\lambda_n), & \alpha < 1 \text{ or } \beta > 0 \\ \left(\frac{12}{1-|a|}\right)\omega_\alpha(1/\lambda_n)(1+\log\lambda_n), & \alpha = 1 \text{ and } \beta = 0 \end{cases}$$

whence (5.2) obviously follows.

**Proof of Theorem 4.** Let $0 \leq \beta < \alpha \leq 1$ and $f \in H^{\omega_\alpha}(\mathbb{T})$. By using (4.3) inequality with taking $A_n = V_n^a, \phi = \omega_\beta(\delta)$ and $\|f\|_\varphi = \|f\|$. Then (4.3) gives the following inequality

$$\|V_n^a(f) - f\|_{\omega_\beta} \lesssim \|V_n^a(f) - f\|_\infty \left(1 + \frac{2}{\omega_\beta\left(\frac{1}{n}\right)}\right) + \sup_{0<\delta\leq\frac{1}{n}} 2\frac{\omega(f,\delta)}{\omega_\beta(\delta)}(1 + \|V_n^a\|).$$

For $\alpha < 1$ from Theorem 1 ve Theorem 3, we have

$$\|V_n^a(f) - f\|_{\omega_\beta} \lesssim \left(\frac{12}{1-|a|}\right)\omega_\alpha(1/\lambda_n)\left(1 + \frac{2}{\omega_\beta(1/n)}\right) + \sup_{0<\delta\leq\frac{1}{n}} 2\frac{\omega_\alpha(\delta)}{\omega_\beta(\delta)}\left(1 + \left(3 + \log\frac{2n-\lambda_n}{\lambda_n}\right)\right).$$

Since $f \in H^{\omega_\alpha}(\mathbb{T})$,

$$\|V_n^a(f) - f\|_{\omega_\beta} \lesssim \left(\frac{12}{1-|a|}\right)\omega_\alpha\left(\frac{1}{\lambda_n}\right)\left(1 + \frac{2}{\omega_\beta\left(\frac{1}{n}\right)}\right) + \sup_{0<\delta\leq\frac{1}{n}} 2\frac{\omega_\alpha\left(\frac{1}{n}\right)}{\omega_\beta\left(\frac{1}{n}\right)}\left(4 + \log\frac{2n-\lambda_n}{\lambda_n}\right).$$

For $\alpha = 1$, we have

$$\|V_n^a(f) - f\|_{\omega_\beta} \lesssim \left(\frac{12}{1-|a|}\right)\omega_\alpha\left(\frac{1}{\lambda_n}\right)(1+\log\lambda_n)\left(1 + \frac{2}{\omega_\beta\left(\frac{1}{n}\right)}\right) + 2\frac{\omega_\alpha\left(\frac{1}{n}\right)}{\omega_\beta\left(\frac{1}{n}\right)}\sup_{0<\delta\leq\frac{1}{n}}\left(4 + \log\frac{2n-\lambda_n}{\lambda_n}\right).$$

Therefore, we have

$$\|f - V_n^a(\lambda, f)\|_{\omega_\beta} \lesssim \begin{cases} \left(\dfrac{12}{1-|a|}\right)\dfrac{\omega_\beta(1/n)}{\omega_\alpha(1/n)}\omega_\alpha(1/\lambda_n)\log\dfrac{2n}{\lambda_n}, & \alpha < 1 \text{ or } \beta > 0 \\ \left(\dfrac{12}{1-|a|}\right)\dfrac{\omega_\beta(1/n)}{\omega_\alpha(1/n)}\omega_\alpha(1/\lambda_n)(1+\log\lambda_n)\log\dfrac{2n}{\lambda_n}, & \alpha = 1 \text{ and } \beta = 0 \end{cases}$$

whence (5.3) obviously follows.

## 7. CONCLUSION

The nonlinear Fourier basis $\{e^{ik\theta_a(t)}, k \in \mathbb{Z}\}$ defined by the nontangential boundary value of the Möbius transformation as a typical family of mono-component signals. This basis has attracted much attention in the field of nonlinear and nonstationary signal processing in recent years. In the present paper, we give rate of convergence of $n$ th generalized de la Vallèe-Poussin mean $V_n^a(\lambda, f)$ of series by nonlinear Fourier basis. Furthermore approximation problems for $n$th generalized de la Vallèe-Poussin mean of series by nonlinear Fourier basis are investigated in the uniform norms and in the generalized Hölder norms.